\documentclass{amsart}
\usepackage{amsfonts}

\begin{document}

\title{ Remark on a conjecture of conformal transformations of Riemannian manifolds }

\begin{center}
\author{A. Raouf  Chouikha  }

\thanks{Universite Paris 13 LAGA UMR 7539, Villetaneuse 93430}
\thanks{chouikha@math.univ-paris13.fr}
 
\begin{abstract}
 Ejiri [E] gave a negative answer to a conjecture of Lichnerowicz concerning Riemannian manifolds with constant scalar curvature admitting an infinitesimal non isometric conformal transformation. With this aim he constructed a warped product of a circle of lenght $T$ and a compact manifold. But he omitted in his analysis the condition that $T$ must to be big enough. Here we give an explicit sharp bound $T_0 < T$ that will make the proof complete. Our presentation is self-contained and mainly uses bifurcation techniques. 
Moreover, we show that there are other such examples and contribute some results to the classification of these manifolds.    
\end{abstract}
\maketitle

\end{center}

\section{Introduction}
Let $(M,g)$ be a $C^\infty$ Riemannian manifold of dimension $n \geq 3.$ Denote by $C(M,g)$ the group of $C^\infty$ conformal diffeomorphisms of $M$ and by $I(M,g)$ the group of isometries of $M$.
It is well known that $C(M,g)$ is a Lie group with respect to the compact open topology.\\ 
Let $C_0(M,g)$ (resp $I_0(M,g)$) denote the connected component of the identity of $C(M,g)$ (resp $I(M,g)$). A subgroup $G$ of $C(M,g)$ is said to be essential if there does not exist a smooth function $\rho$ for which $G \subset I(M,e^{2\rho}g)$.\\
A. Lichnerowicz asked the following question : {\it are there manifolds $(M,g)$ which are not conformal to Euclidean space or to the standard sphere $(S^n,g_0)$ for which $C_0(M,g)$ is essential} ?\\
In the compact case J. Lelong-Ferrand [LF] and M. Obata [O] proved the following \\

{\bf Proposition 1:} \quad {\it Let $(M,g)$ be a $C^\infty$ compact Riemannian manifold of dimension $n \geq 3$. If the subgroup $C_0(M,g)$ is essential and the scalar curvature of $(M,g)$ is constant then $(M,g)$ is isometric to the standard sphere $(S^n, g_0)$.}\\

However, one espects that condition "$C_0(M,g)$ is essential" may be replaced by a more generalized condition "$C_0(M,g) \neq I_0(M,g)$". This means:  if $(M,g)$ admits a infinitesimal non isometric conformal transformation then $(M,g)$ is isometric to the standard sphere.\\
Various improvements were brought, notably when $(M,g)$ is a locally conformally flat manifold rather than of constant scalar curvature, or when the Ricci tensor of $(M,g)$ is parallel. Nevertheless, it was necessary to await the counter example of Ejiri [E] to exhaust this research direction. 

\section{A counter example}

Let us consider \ $X = f(\frac {d}{dt})$\  an infinitesimal conformal transformation of $(S^1\times _f N,\tilde g)$ verifying ${\it L}_X \tilde g = 2 f' \tilde g$ where $\frac {d}{dt}$ is a vector field on $S^1$, ${\it L}_X $ is the Lie differentiation in the direction of $X$. \\
Let $(N_1,h_1)$  and $(N_2,h_2)$ be Riemannian manifolds, $f : N_1 \rightarrow R$ a positive function. Define the $f$-warped product $N_1\times _f N_2$ of $N_1$ and $N_2$ to be the Riemannian manifold $N_1\times N_2, h_1+fh_2$ with $$(h_1+fh_2)_{(x,y)}(u+X,v+Y) = h_1(u,v) + f(x)h_2(X,Y)$$
for $u,v \in T_x N_1; X,Y \in T_y N_2.$\\

Let $S^1$ be the circle parametrized by arc-lenght $t$ of lenght 
$$T = \int_{S^1} dt.$$
This circle is endowed with the Euclidean metric $dt^2$ and let $f$ be a $C^\infty$ function on $S^1$. \\
Consider $(N,h)$ be a (n-1)-dimensional compact Riemannian manifold with positive constant scalar curvature \ $R$ .\  Let \ $\tilde R$ \ be the scalar curvature of the warped product\ $(S^1\times _f N, \tilde g)$\ where\ $\tilde g = dt^2 + f(t)h$. \ Then the scalar curvatures and\ $f$\ satisfy the equation
$$ (E)\qquad  \tilde R f^2 + 2(n-1) f f'' + (n-1)(n-2) f'^2 - R = 0$$
At first remark that Equation (E) always admits 
$$f = constant$$ as a trivial periodic solution.

Notice that \ $T = \int_{S^1} dt$\ must precisely be the period of the (non constant) solution \ $f$.\  
The corresponding manifolds are bundles with fibres \ $N$
\ over the circle \
$S^1$ \ (parametrized by arc length \ $t$ ) \ equipped with the warped metrics.\\
Ejiri remarked that if \ $\tilde R$ \ is a positive constant and \ $f$ \ is a positive periodic solution of (E), then \ $(S^1\times _f N, \tilde g)$ \ is a compact Riemannian manifold of dimension $n$ admitting an infinitesimal non isometric conformal transformation.\\

Notice that the solution $f$ must be non constant otherwise \ ${\it L}_X \tilde g = 2 f' \tilde g = 0$\ and the manifold \ $(S^1\times _f N, \tilde g)$ \ does not admit an infinitesimal non isometric conformal transformation.\\

Thus,  Ejiri constructs a warped product of a circle and a compact Riemannian manifold with constant scalar curvature, [E]. More exactly, he asserts the following\\

{\bf Proposition 2} \quad {\it Let $(N,h)$ be a compact Riemannian manifold of dimension $n-1$ with positive scalar curvature. Then there exists a positive periodic function $f$ on a circle $(S^1, dt^2)$ such that the warped product $(S^1\times N,dt^2+f(t)h)$ has  a constant scalar curvature and admits an infinitesimal non isometric conformal transformation.}\\ 

But another condition is necessary to ensure that the solution $f$ is periodic non constant.

\newpage

{\bf Analysis of Equation (E)}\\

Now consider the change $$f = \alpha (1+ j)^{2/n}$$ where \ $\alpha = (\frac {R}{\tilde R})^{\frac {n}{4}}.$ Then $$f' =  \frac {2\alpha j'}{n} (1+j)^{\frac {2}{n}-1}$$ and $$ f'' = (\frac {2}{n}-1)\frac {2\alpha }{n} (1+j)^{\frac {2}{n}-2} j'^2 + \frac {2}{n}\alpha (1+j)^{\frac {2}{n}-1} j''.$$
We calculate  $$2 f f'' + (n-2) f'^2 = 2 \alpha (1+ j)^{2/n} [(\frac {2}{n}-1)\frac {2\alpha }{n}(1+j)^{\frac {2}{n}-2} j'^2 + \frac {2}{n}\alpha (1+j)^{\frac {2}{n}-1} j''] +$$ $$ \frac {2(n-2)\alpha ^2 j'^2}{n} (1+j)^{\frac {4}{n}-2}.$$
After simplication one finds
$$2 f f'' + (n-2) f'^2 = \frac {4}{n}\alpha^{\frac {4}{n}-1} \alpha^{\frac {n}{2}} (1+j)^{\frac {4}{n}-1} j''$$
Replace in Equation (E) which becomes

\begin{equation}
 j'' -  \frac {n \tilde R}{4(n-1)} (1+j)^{1-4/n} = - \frac {n \tilde R}{ 4(n-1)} (1+j)
\end{equation}

But, as we will prove below, a non constant periodic function does not exist for any circle lenght\ $T$.\ For example, if \ $T$\ takes a small value then Equation $(1)$ admits only constants as periodic solutions. 
More precisely, when \ $j$ \ is closed to \ $0$\ this equation can be written 
$$ (1') \qquad  j'' -  \frac {n \tilde R}{4(n-1)} (1+(1-4/n)j +(-2/n) j^2 +....) +{n\tilde R\over 4(n-1)} (1+j) = 0$$
where the constant $\alpha =  \bigg({R\over{ {\tilde R}}}\bigg)^{n/4}.$\\
Then \ $j$\ is closed to the solution of the linearized equation  

$$ (1') \qquad j'' + \frac {\tilde R}{n-1} j = 0$$
 
This means Equation (1) bifurcates at $j \equiv 0$ when $\frac {\tilde R}{n-1} = (\frac {2 \pi}{T})^2$, \ [C-R]\\ 

Thus, there is a positive bound \ $T_0$\  such that if \ $T \leq T_0$\ the above equation may have only constant solutions, i.e.
\ $f(t) \equiv \alpha = \bigg({R\over{4\tilde R}}\bigg)^{n/4}.$\\ 
This means that there is a bound \ $T_0$\ such that condition \ $T > T_0$\ appears to be necessary for a non constant  $T$-periodic solution of $(E)$ exists.\\ Hence, examples given by Ejiri lack precision and his proof is incomplete.

\newpage

\section{Condition on the scalar curvature of $(S^1\times _f N, \tilde g)$}

In this section we shall add another necessary condition so that Proposition 2 becomes correct.
In particular, the existence of  non constant periodic solutions $f$ depend on the scalar curvature $R$ and the lengh $T$ of the circle. The following result gives a precise bound $T_0$ such that if \ $T \leq T_0$\ the function $f$ cannot be (non constant) periodic solution of $(E)$\\

{\bf Theorem 3} \quad {\it Let $(N,h)$ be a compact Riemannian manifold with positive scalar curvature $R$ of dimension $(n-1), n \geq 3$. Let $S^1$ be the circle of lenght $T$ and $\tilde R$ be a positive constant verifying the condition 
$$T > \frac {2\pi \sqrt {n-1} }{\sqrt {\tilde R}} = T_0.$$
Then the differential equation 
$$ \tilde R f^2 + 2(n-1) f f'' + (n-1)(n-2) f'^2 - R = 0$$
admits a positive non constant periodic solution\ $f$.\\
Moreover, the warped product $(S^1\times _f N, dt^2 + f(t)h)$ is a compact Riemannian manifold with scalar curvature $\tilde R$ of dimension $n$ admitting an infinitesimal non isometric conformal transformation.}
\bigskip

\section{Proof}
Now we again make the change of function $$f = x^{2/n}$$ then $$f' = \frac {2}{n} x^{\frac {2}{n}-1} x'\quad and \quad f'' = (\frac {2}{n}-1)\frac {2}{n} x^{\frac {2}{n}-2} x'^2 + \frac {2}{n} x^{\frac {2}{n}-1} x''.$$
We calculate  $$2 f f'' + (n-2) f'^2 = (\frac {2}{n}-1)\frac {4}{n} x^{\frac {4}{n}-2} + (n-2)\frac {4}{n^2} x^{\frac {4}{n}-2} + \frac {4}{n} x^{\frac {4}{n}-1} x'' = \frac {4}{n} x^{\frac {4}{n}-1} x''$$
Replace in the following equation  

\begin{equation} 
\tilde R f^2 + 2(n-1) f f'' + (n-1)(n-2) f'^2 - R = 0
\end{equation}

which becomes

\begin{equation} 
 x'' - \frac {nR}{4(n-1)}x^{1-4/n} = -\frac {n\tilde R}{ 4(n-1)} x
\end{equation}

The last equation has been analysed in studying the parallelism of the Ricci tensor of the manifold $S^1\times _f N, \tilde g$ where $N$ is a Einstein manifold, see [Ch] and [D].\\

 Let \ $c$\ be the energy level for that equation. All periodic orbits \ $\gamma_c(t)$ \ of the following system which is equivalent to Equation (3) 

\begin{equation}   
\begin{cases}x' = - y  & \cr 
y' =  \frac {nR}{4(n-1)} x^{1-\frac {4}{n}} - \frac {n\tilde R}{ 4(n-1)} x,  & \cr
\end{cases}  
\end{equation}

are surrounded by the homoclinic orbit \
$\gamma_{c_0}$ .\ The latter one may be parametrized by \ $(x_0(t),y_0(t)).\ x_0(t)$\ is the (degenerate) non periodic solution of (3).\\
Denote the coordinates of \ $\gamma_c(t)$ \ by \ $(x_c(t),y_c(t).$ \
 When the value \ $c$ \ satisfies the condition \ $1 < c < c_0$ \ , the
correspondant orbit is periodic ($c_0$\ corresponds to a periodic solution of null energy) .\\
The center of System (4) is \ $(x=\alpha ,y=0)$, \ where \ $\alpha  =  (\frac {R}{4\tilde R})^{4/n}$.\\
 One may easily remark that two positive T-periodic solutions of
(4) having the same energy translate of one another, and thus give rise to equivalent
metrics on \
$(S^1(T)\times N),g_0)$. \\  Note that the metric corresponding to the conformal factor
\quad
$u_0\ : \quad  \ g = {x_0}^{\frac {4}{n-2}}g_0$ \quad is non complete. Therefore the constant \ $c$ \ cannot attain the critical value \ $c_0$ \ .\\

Equation (3) may be written under the following form

\begin{equation}
x'' + \phi (x) = 0
\end{equation}
where
$$\phi (x) = \frac {nR}{{4(n-1)}} (x - \alpha ) - \frac {nR}{{4(n-1)}}(x - \alpha )^{1-4/n} .$$
The period of the periodic solutions depends on the energy \ $T\equiv T(c)$\  with
$c$ the energy constant. It can be expressed as
$$T(c) = {\sqrt 2} \int_a^b \frac {du }{\sqrt{c - G(u)}}$$
where $G(u)$ is an integral of $\phi (u),$ with a nondegenerate relative minimum
at the origin. It satisfies in addition,  \ $ G(a) = G(b) = c$\ and \ $ a \leq \alpha \leq b$.  \\
So, \ $\phi (\alpha ) = 0 $\ and \ $\phi '(\alpha ) = \frac {nR}{{4(n-1)}} > 0.$\\

 Hence, \ $x = \alpha $\  is a center for Equation (3)i.e. in the neighbourhood of the trivial solution \ $h(t) \equiv \alpha$\ Equation (2) admits  a periodic solution.\\ 
The following lemma is a classical result of global bifurcation theory (for details see for example [C-R]).

\bigskip

 {\bf Lemma 4}\qquad {\it Under the above hypothesis the  family of solutions \ $(T,u_T(t))$ \ of the ODE (3) 
(where \ $T$ \ is the minimal period) has bifurcation points on the values
\quad $(T_k,u_T{_k}(t))$ \ where \ $ T_k=\frac {{2\pi k}}{{\sqrt{n-2}}}$\ and \ $u_{T_k}\equiv\alpha$\ is a constant.  
In this family, there is a curve of non trivial solutions which bifurcates to the right of the
trivial one.}\\

\smallskip
So, let us consider a positive T-periodic solution: if \ $T\neq T_k$\ , then the linearized associate
equation is non-singular.\\ We may also
deduce from the bifurcation theorem, applied to the simple eigenvalues problem, that there is a
unique curve of non trivial solutions near the point \ $(T_k,\alpha).$ \  By Lemma 4, this uniqueness is
global. The trivial curve is \ $u_T \equiv\alpha$ .\\ Moreover, \
${\displaystyle (\frac {du}{{dT}})_{T=T_k}}$ \ is an eigenvalue of the corresponding linearized equation.
According to global bifurcation theory, we assert that the non trivial curves turn on the
right of the singular solution \ $(T_k,u_{T_k}(t))$ \ . Consequently, when \ $T$ \ varies, two non
trivial curves never cross.\\ Thus, in particular for \ $T\leq T_0$\ Equation (3) does not admit a non trivial
 periodic solution.
 
\subsection{General remarks}
One can refine the study of Equation $(E),$ for example by looking for the exact number of periodic solutions for that equation. This is more interesting of the ODE point of view,  but not so much regarding the geometric problem presented at first.  What is important here is to notice that when the following condition on T holds. 
$$T \leq \frac {2\pi \sqrt {n-1} }{\sqrt {\tilde R}}$$
then the warped product $(S^1\times _f N, dt^2 + f(t)h)$ does not admit any infinitesimal non isometric conformal transformation (because Equation (E) has only constants as periodic solutions).\\
We have already examined the existence and the number of these metrics [Ch].

\bigskip

\section {Other examples}

An interesting question is the following : \\ 

{\it Are there other manifolds non isometric to the standard sphere admitting an infinitesimal non isometric conformal transformation ?}\\

We are able to produce a partial answer using a result of Derdzinski [D] concerning the classification of all n-compact Riemannian manifolds \ $(M,g), \ n\geq 3$\ with harmonic curvature.\\
Notice that such a manifold displays other interesting geometrical properties.
In the case where $(N,g_0)$ is an Einstein manifold.
Derdzinski [D] established a classification of the compact n-dimensional Riemannian
manifolds\quad $(M_n,g), \ n \geq 3$, \quad with harmonic curvature. If the Ricci tensor \ $Ric(g)$ \
is not parallel and has less than three distinct eigenvalues at each point, then \ $(M,g)$ \ is isometrically 
covered  by a manifold  $$(S^1(T) \times N,dt^2 + h^{4/n}(t)g_0),$$ where
the non constant positive periodic solutions \ $h$ \ satisfies Equation (E). Here \ $(N,g_0)$\ is a (n-1)-
dimensional Einstein manifold with positive (constant) scalar curvature.\\

We thus obtain the following

\bigskip

{\bf Theorem 4}\quad {\it Let \ $(M,g)$\ a compact Riemannian manifold of dimension \ $n, \\ n \geq 3$ \ with harmonic curvature and covered isometrically by a manifold  $$(S^1(T) \times N,dt^2 + h^{4/n}(t)g_0),$$ where
the non constant positive periodic solutions \ $h$ \ satisfies Equation (E). Suppose that its Ricci tensor is not parallel and has less than three distinct eigenvalues at each point; then, \ $(M,g)$\ admits an infinitesimal non isometric conformal transformation.}

\newpage

{\it Acknowledgements} : I would like to thank the referee for helpful remarks. 

\bigskip
{\bf References}\\

[E] N. Ejiri {\it A negative answer to a conjecture of conformal transformations of Riemannian manifolds} J. Math. Soc. Japan, vol 33, No 2, p.261-265, 1981.\\

[O]M. Obata {\it The conjecture of conformal transformations of Riemannian manifolds} J. diff Geom., vol 6, p. 247-258, 1971.\\

[LF]Lelong-Ferrand {\it Transformations conformes et quasi-conformes des variétés
riemanniennes compactes} Acad. Roy. Belg., Mém. Coll. 39, no. 5,1971.\\

[ch]A. R. Chouikha {\it On the existence of metrics with harmonic curvature and non parallel Ricci tensor}  Balk. J. of Diff. Geom. and Appl., p. 21-31, vol 8, 2003\\

[D]A. Derdzinski {\it On compact Riemannian manifolds with harmonic curvature }  Math.
Ann., vol 259, (1982), 145-152.\\

[C-R] M.Crandall and P.Rabinowitz \ {\it Bifurcation from simple eigenvalues } J. Funct.
Anal. vol 8, (1971), 321-340.\\

\end{document}